\documentclass[11 pt]{amsart} 

\usepackage{graphicx}
\usepackage{amsmath,amsthm,amscd,amssymb}
\usepackage{amsfonts}
\usepackage{color}
\usepackage{latexsym}
\usepackage{epsfig,epstopdf}
\usepackage{enumerate}

\newcommand{\ds}{\displaystyle}

\newcommand{\R}{\mathbb{R}}            
\newcommand{\C}{\mathcal{C}}	

\newcommand{\Z}{\mathcal{Z}}
\newcommand{\A}{\mathcal{A}}
\newcommand{\X}{\mathcal{X}}
\newcommand{\U}{\mathcal{U}}

\newcommand\norm[1]{\left\| #1\right\|}
\newcommand\inner[1]{\langle #1\rangle}

\newcommand{\Aa}{\mathbb{A}}
\newcommand{\Ba}{\mathbb{B}}
\newcommand{\Ia}{\mathbb{I}}
\newcommand{\Ga}{\mathbb{M}_g}
\newcommand{\Fa}{\mathbb{F}}

\newtheorem{definition}{\sc Definition}[section]
\newtheorem{teo}{\sc Theorem}[section]
\newtheorem{prop}{\sc Proposition}[section]

\newtheorem{lemma}{\sc Lemma}[section]

\newcommand{\dem}{\noindent {\bf Proof.} \mbox{}}

\begin{document}

\title[Controllability of the Semilinear Beam Equation]
{Controllability of the Impulsive Semilinear Beam Equation with Memory and Delay}
\date{\today}
\author[A. CARRASCO, C. GUEVARA  AND H. LEIVA ]{A. Carrasco$^1$, C. Guevara$^2$ and H. Leiva$^3$  }
\address{$^{1}$ Universidad Centroccidental Lisandro Alvarado \\
          Decanato de Ciencias y Tecnolog\'ia, Departamento de Matem\'aticas \\
          Barquisimeto 3001-Venezuela}
          \email {acarrasco@ucla.edu.ve}
\address{$^{2}$ Louisiana State University \\
          College of Science, Department of Mathematics \\
          Baton Rouge, LA 70803-USA} \email{cguevara@lsu.edu, cristi.guevara@asu.edu}
\address{$^{3}$ School of Mathematical Sciences and Information Technology\\
           Universidad Yachay Tech \\
         San Miguel de Urcuqui, Ecuador}
          \email{hleiva@yachaytech.edu.ec, hleiva@ula.ve}


\begin{abstract}
The semilinear beam equation with impulses, memory and delay is considered and its approximate controllability obtained. This is done by employing a technique that avoids fixed point theorems and pulling back the control solution to a fixed curve in a short time interval. Demonstrating, once again, that  the controllability of a system is robust under the influence of impulses and delays.

\end{abstract}

\maketitle 


\section{Introduction}

Beams have been used since ancient times to reinforce structures such as bridges, buildings, and others. Through the millennia, understanding the dynamics and controllability of beams, including bending and vibration has been  of great importance. Pioneering studies goes back to 1493 Leonardo da Vinci's manucript that identified properly the stresses and strains in a beam subject to bending \cite{DaVinci} and Galileo Galilei's writings that identified the principle of virtual work as a general law but made incorrect assumptions \cite{Timoshenko:1921aa}. It was not until the late 17th century with the elasticity theory evolution that  Leonhard Euler and Daniel Bernoulli provided  a second-order spatial derivatives mathematical model that later, in 1921, Stephen Timoshenko improved by including a shear deformation and rotational inertia effects, obtaining  fourth order  mathematical model (see \cite{Timoshenko:1921aa,Timoshenko:1922aa, Timoshenko:1953aa} for  details).

Nowadays, adjustments of the Timoshenko beam model, in  mechanical engineering and nanotechnology design \cite{Wang-Tan:2006aa,Wang-Zhang:2006aa}, yield to the impulsive semilinear beam equations of the form \eqref{eq:beam} where  the memory and delay provide information of the viscoelasticity property  and response of the materials. 

In this paper, we are exploring the  approximate  controllability on a bounded domain $\Omega \subseteq \R^{N} \,(N\geq1)$ of
\begin{equation}\label{eq:beam}
 w_{tt}  -  2\beta\Delta w_t + \Delta^{2}w = u(t,x) + f(t,w(t-r),w_{t}(t-r),u)+\ds \int_{0}^{t}M(t-s)g(w(s-r,x))ds,
\end{equation}
subjected to the initial-boundary conditions and impulses
\begin{equation}\label{eq:initial}
\left\{\begin{array}{ll}
 w(t,x)=\Delta w (t,x)=0, &\mbox{in}\; (0, \tau) \times \partial \Omega,\\
 \begin{split}
 &w(s,x)=\phi_1(s,x),\\
 &w_{t}(s,x)=\phi_2(s,x),
 \end{split} & \mbox{in} \; [-r,0] \times \Omega,\\
 w_{t}(t_{k}^{+},x) = w_{t}(t_{k}^{-},x)+I_{k}(t_k,w(t_{k},x),w_{t}(t_{k},x),u(t_{k},x)), &t \neq t_k, \; k=1, \dots, p,
  \end{array}
 \right.
\end{equation}
where $\Delta w=\sum_{j=1}^{N}\frac{\partial^2w}{\partial x_j^2}$ and $\Delta^2 w=\sum_{j=1}^{N}\frac{\partial^4w}{\partial x_j^4}$. Additionally, the  damping coefficient  $\beta > 1$  and the real-valued functions $w = w(t,x)$  in $(0, \tau] \times  \Omega$ represents the beam deflection, $u$ in $(0, \tau] \times  \Omega$ is the distributed control, $M$ acts as convolution kernel with respect to the time variable, the impulses $I_k$ are defined on $[0, \tau] \times \R^3$  and  the nonlinearities $g$ on $\R$,  $f$ on $[0, \tau] \times \R^3$. Under the assumptions:
\begin{description}
\item[H$_1$] $M\in L^{\infty}((0,\tau)\times \Omega)$,  and $g, f, I_{k}$  are smooth enough, in order that, for all $\phi, \psi\in \C([-r,0],L^{2}(\Omega))$ and $u\in L^{2}([0,\tau]; L^{2}(\Omega))$ the equation \eqref{eq:beam} admits only one mild solution on $[-r,\tau]$.
\item [H$_2$] $t\!\in\! [0,\tau],$\; $a, b \geq \!0$ and  $u, v,  y \in \R$, the nonlinearity $f$ satisfies
\begin{equation}\label{eq:f1}
    \begin{array}{ll}
|f(t,y,v,u)|  & \leq a\sqrt{|y |^2+ |v |^2} +b.
\end{array}
\end{equation}
\end{description}

N.~Abada, M.~Benchohra, and H.~Hammouche \cite{Abada-Benchohra:2010aa} and  R.~S. Jain and M.~B. Dhakne in \cite{Jain-Dhakne:2013aa} works showed the existence of solutions for impulsive evolution equations with delays. Balachandran,  Kiruthika, and Trujillo \cite{Balachandran-Kiruthika:2011aa} supplied  existence results for the fractional impulsive integrodifferential equations and finally for  the Beam equation with variable coefficients  J.~L{\'\i}maco, H.~Clark, and A.~Feitosa \cite{Limaco-Clark:2005aa} showed   the existence and uniqueness of non-local strong solutions and the  existence of a unique global weak solution with decay rate energy.  

Inspired in a series of papers from A. Carrasco, H. Leiva, N. Merentes and J.
Sanchez on the approximate controllability of semilinear beam equation  \cite{Carrasco-Leiva:2013aa,Carrasco-Leiva:2014aa,Carrasco-Leiva:2016aa} and the works on the approximate controllability for the semilinear heat and strongly damped wave equations with memory and delays by C. Guevara and H. Leiva  \cite{Guevara-Leiva:2016aa, Guevara-Leiva:2017aa}. We prove the  approximate controllability of the beam equation  \eqref{eq:beam} under the initial-boundary condition \eqref{eq:initial}  with memory, impulses and delay terms  by applying   A.E. Bashirov, N. Ghahramanlou, N. Mahmudov, N. Semi and H. Etikan technique \cite{Bashirov-Ghahramanlou:2014aa, Bashirov-Ghahramanlou:2015aa, Bashirov-Jneid:2013aa, Bashirov-Etikan:2010aa},  and avoiding the Rothe's fixed point theorem used in \cite{Carrasco-Leiva:2013aa,Carrasco-Leiva:2016aa} and the   Schauder fixed point theorem applied in \cite{Carrasco-Leiva:2014aa}.

The structure of this paper is as follow: In section \ref{sec:formulation}, we present the abstract formulation of the beam equation \eqref{eq:beam}. Section \ref{sec:lineal}, recalls the linear controllability characterization of the problem. In section \ref{sec:semilineal}, the approximated controllability of the beam equation with memory, delay and impulses is proved.

\section{Abstract Formulation of the Problem} \label{sec:formulation}

In this section, we choose the appropiate Hilbert space where the Cauchy problem \eqref{eq:beam}-\eqref{eq:initial}  can be written as an abstract differential equation.

 First of all, notice that the term $-2\beta\Delta w_t$ in the equation \eqref{eq:beam} acts as a damping force, thus the energy space used to set up the wave equation is not suitable here. Even so,  in   \cite{Oliveira:1998aa}, Oliveira shows that the uncontrolled linear  equation can be transformed into a system of parabolic equations of the form $w_{t} = D \Delta w$,  obtaining that corresponding space for the abstract formulation of the problem is $\Z^{1}=\left[H^{2}(\Omega) \bigcap H^{1}_{0}(\Omega) \right] \times L^{2}(\Omega)$ and proving that the linear part of this system generates a strongly continuous analytic semigroup in this space.

Consider the Hilbert space $\X = L^{2}(\Omega)$,  and denote  $\A=-\Delta$ with eigenvalues  $0 <
\lambda_{1}<\lambda_{2}<...<\lambda_{j}\to \infty,$ with multiplicity $\gamma_{j}<\infty$ equal to its corresponding eigenspace dimension. Recall, $\A$ satisfies the following properties:
\begin{enumerate}[(i)]
\item There exists a complete orthonormal set $\left\{
\phi_{j_k} \right\}$ of eigenvectors of $\A$.

\item For all $x \in D(\A)$,
\begin{equation*} \label{prop}
\A x = \sum_{j = 1}^{\infty} \lambda_{j} \sum_{k =
1}^{\gamma_j} \inner{\xi, \phi_{j_k}} \phi_{j_k} =\sum_{j = 1}^{\infty}
\lambda_{j}E_{j}x,
\end{equation*}
where $\inner{\cdot, \cdot}$ denotes the inner product in $\X$, $\ds E_{n}x = \sum_{k = 1}^{\gamma_j} \inner{z, \phi_{j_k}} \phi_{j_k},$ and  $\{ E_j \}$ is a family of complete orthogonal projections in
$\X$.
 \item $-\A$ generates an analytic semigroup $\{
S(t) \}_{t \geq 0}$ given by
$$
S(t)x =  \sum_{j = 1}^{\infty} e^{-\lambda_j t}E_{j}x \quad  \mbox{and} \quad  \norm {S(t)} \leq e^{-\lambda_{1}t}.
$$
\item For $\alpha\geq0$ the fractional powered spaces $\X^{\alpha}$ are given by
\begin{equation*}
   \X^{\alpha} =D(\A^{\alpha}) = \left\{x \in \X : \sum_{j = 1}^{\infty} \lambda_{j} ^{2
\alpha} \norm{ E_{j}x}^2 < \infty \right\}
\end{equation*}
equipped with the norm $\ds
\norm{x}_{\alpha}^2 = \norm{\A^{\alpha}x}^2=  \sum_{j = 1}^{\infty}
\lambda_{j}^{2 \alpha} \norm{ E_{j}x}^2
$, where $\ds \A^{\alpha}x = \sum_{j = 1}^{\infty} \lambda_{j}^{ \alpha}  E_{j}x$.
\end{enumerate}
In particular, $\alpha=2$ yields $\ds \A^{2}x = \sum_{j = 1}^{\infty} \lambda_{j}^{ 2}  E_{j}x=(-\Delta)^2x=\Delta^2x$. And for $\alpha=1,$ the Hilbert space $\Z^{1}=\X^{1}\times \X$ has the  norm
$$
\ds \norm {\left(
                  \begin{array}{c}
                    w \\
                    v \\
                  \end{array}
                \right)}^{2}_{\Z^{1}}=\|w\|^{2}_{1}+\|v\|^{2}.
$$

Using the above notation, we rewrite the system \eqref{eq:beam}-\eqref{eq:initial}  as the second-order ordinary differential equations in the Hilbert
space $\X$
\begin{equation} \label{eq:ODE}
\left\{%
\begin{array}{ll}
    \begin{split}
    w''(t) = &-\A^{2}w(t) - 2\beta \A w'(t) + u(t)+  \ds \int_{0}^{t}M(t,s)g^{e}(w(s-r))ds  \\
 &\ \ + f^{e}(t,w(t-r),w'(t-r),u(t)),
 \end{split}
  & t>0,t \neq t_k,\\
  \begin{split}
     w(s)  &=  \phi_1(s), \\
      w'(s)&=\phi_2(s),
      \end{split}
 &s \in [-r,0],\\
 w' (t_{k}^{+}) = w' (t_{k}^{-})+I^{e}_{k}(t_k,w(t_{k}),w' (t_{k}),u(t_{k},)), & k=1, \dots, p,
    \end{array}%
\right.
\end{equation}
where $\U=\X=L^{2}(\Omega)$, and
\begin{eqnarray*}
I_{k}^{e}:&[0, \tau]\times \Z^{1} \times \U &\longrightarrow \qquad \X \\
&(t,w,v,u)(\cdot)&\longmapsto \quad I_{k}(t,w(\cdot),v(\cdot),u(\cdot)),
\end{eqnarray*}
\begin{eqnarray*}
f^{e}:&[0, \tau]\times \C (-r,0;  \Z^{1} ) \times \U & \longrightarrow \qquad \X\\
&(t,\Phi,u)(\cdot)&\longmapsto \quad f(t,\phi_1(-r, \cdot),\phi_2(-r, \cdot),u(\cdot)),
\end{eqnarray*}
and
\begin{eqnarray*}
g^{e}:&\C(-r,0;  \Z^{1} )  &\longrightarrow  \Z^{1} \\
&\Phi=\left(\begin{array}{c}
             \phi_1\\
             \phi_2
        \end{array}\right)&\longmapsto  g(\phi_1(\cdot-r)).
\end{eqnarray*}

\noindent Changing variables, $v=w',$  the systems \eqref{eq:ODE} can be written as an abstract first order
functional differential equations with memory, impulses and delay in $\Z^{1}$
\begin{equation}\label{eq:abstract}
\left\{
\begin{array}{lr}
z' =  -\Aa z+ \Ba u + \ds \int_{0}^{t}\Ga(t,s ,z_{s}(-r))ds + \Fa(t,z_{t}(-r),u(s)) ,& z\in Z^{1},\; t\geq 0, \\
z(s)  = \Phi(s), & s \in [-r,0],\\
z(t_{k}^{+})  = z(t_{k}^{-})+\Ia_{k}(t_k, z(t_{k}),u(t_{k})), &  k=1,2,3, \dots, p,
\end{array}
\right.
\end{equation}
where $ z =\left(\begin{array}{c}
             w\\
             v
        \end{array}\right)$,
   $\Phi=\left(\begin{array}{c}
             \phi_1\\
             \phi_2
        \end{array}\right) \in \C\left(-r,0;  \Z^{1} \right),$ $u\in L^{2}(0,\tau;\U)$, $\Aa = \left(
   \begin{array}{rr}
     0 & I_{\X} \\ - \A^2 & -2\beta \A
   \end{array}\right)$ is a unbounded linear operator  with domain
$$
D(\Aa)=\{w\in H^{4}(\Omega):\:w=\Delta w=0\}\times D(\A),
$$ and  $I_{\X}$ being the identity in $\X $.
   $\Ba: \U \longrightarrow \Z^{1}$ is the bounded linear operator defined by
$\Ba u=
\left(\begin{array}{c}
             0\\
             u
        \end{array}\right),$ and the functions
\begin{eqnarray*}
\Ia_{k}:&[0, \tau]\times \Z^{1} \times \U& \longrightarrow \qquad \Z^{1} \\
&(t, z,u)&\longmapsto \quad
\left(\begin{array}{c}
             0\\
             I_{k}^{e}(t,w,v,u)
        \end{array}\right)
        \end{eqnarray*}
\begin{eqnarray}
\label{functionF}
\Fa:& [0, \tau] \times \C(-r,0;  \Z^{1} ) \times \U & \longrightarrow \qquad \Z^{1}\\
&(t, \Phi,u)&\longmapsto \quad
\left(
  \begin{array}{c}
    0 \\f^{e}(t,\phi_1(-r),\phi_2(-r),u)
  \end{array}
  \right),\notag
  \end{eqnarray}
and
\begin{eqnarray*}
\Ga:&[0,\tau]\times [0,\tau]\times \C(-r,0;  \Z^{1} )  &\longrightarrow  \Z^{1} \\
&(t,s,\Phi)&\longmapsto \left(
  \begin{array}{c}
    0 \\M(t,s) g^{e}(\Phi)
  \end{array}
  \right).
\end{eqnarray*}
Moreover, this abstract formulation together with condition \eqref{eq:f1} and the continous imbeding $\X^1 \subset \X$ yields
\begin{prop}\label{prop:cotaF}
There exist constants  $\tilde{a},\tilde{b}>0$ such that, for all $(t, \Phi,u) \in [0, \tau] \times \C(-r,0;  \Z^{1} ) \times \U$ the following inequality holds
\begin{equation}\label{eq:bound}
\norm{\Fa(t,\Phi,u)}_{\Z^{1}}  \leq   \tilde{a}\| \Phi(-r) \|_{\Z^1}+\tilde{b}.
\end{equation}
\end{prop}
A.~ Carrasco, H.~ Leiva, and J.~Sanchez \cite[Theorem 2.1]{Carrasco-Leiva:2013aa} proved that the linear unbounded operator
$\Aa$  generates a strongly continuous compact semigroup
 in the space $\Z^1$ which decays exponentially to zero, precisely:
\begin{prop}\label{semigroup}
The operator $\mathcal{A}$ is the infinitesimal generator
of a strongly continuous compact semigroup $\{T(t)\}_{t\geq0}$
represented by
\begin{equation}\label{repre}
T(t)z=\displaystyle\sum_{j=1}^{\infty}e^{\mathbb{A}_{j}t}P_{j}z,\qquad z\in
\Z^{1},\;t\geq 0,
\end{equation}
where $\{P_{j}\}_{j\geq0}$ is a complete family of orthogonal
projections in the Hilbert space $\Z^{1}$ given by
\begin{equation}\label{proyecciones}
    P_{j} = diag(E_{j},E_{j}),
\end{equation}
and
$$
\mathbb{A}_{j}=K_{j}P_{j},\qquad
K_{j}=\left(
                           \begin{array}{cc}
                             0 & 1 \\
                             -\lambda_{j}^{2}  & -2\beta\lambda_{j} \\
                           \end{array}
                         \right),\qquad j\geq 1,
$$
and there exists $M \geq 1$ and $\mu >0$ such that
$$
\parallel T(t)\parallel\leq Me^{-\mu t},\qquad t\geq0.
$$

\end{prop}

\section{Approximate Controllability of the Linear System }\label{sec:lineal}
 This section is devoted to characterize the approximate controllability of the linear system. Thus,  for all $z_{0}\in \Z^{1}$ and $u\in
L^{2}([0,\tau];\U)$ consider the initial value problem
\begin{equation}\label{eq:linear}
\left\{%
\begin{array}{lll}
    z'(t) = \mathbb{A} z(t) + \Ba u(t),\\
    z(t_{0}) = z_{0},
    \end{array}%
\right.
\end{equation}
obtained from \eqref{eq:abstract}. It admits only one mild solution on $0\leq t_{0}\leq t\leq \tau$ given by
\begin{equation}\label{eq:mild-linear}
    z(t)=T(t-t_{0})z_{0} + \displaystyle\int_{t_{0}}^{t}T(t-s)\Ba u(s)ds.
\end{equation}
\begin{definition} \label{def2}
({\bf Approximate Controllability of (\ref{eq:linear})}) The system (\ref{eq:linear}) is said
to be approximately controllable on $[t_{0},\tau]$ if for every $z_0$,
$z_1\in \Z$, $\varepsilon>0$ there exists $u\in L^{2}(t_{0},\tau;\U)$ such
that the solution $z(t)$ of (\ref{eq:mild-linear}) corresponding to $u$
verifies: $$\|z(\tau)-z_1\|<\varepsilon.$$
\end{definition}

For the system \eqref{eq:linear} and $\tau>0$, we have the following notions:
\begin{enumerate}
\item $G_{\tau\delta}$ is the controllability  operator defined by
\begin{eqnarray*}
G_{\tau\delta}: L^2(\tau-\delta,\tau;\U) \longrightarrow& \Z^{1}\\
u\longmapsto&\ds \int_{\tau-\delta}^{\tau}T(\tau-s)\Ba u(s)ds,
\end{eqnarray*}
with corresponding adjoint $G^*_{\tau\delta}$ given by
\begin{eqnarray*}
G^*_{\tau\delta}:  \Z^{1} \longrightarrow& L^2(\tau-\delta,\tau;\U)\\
z\longmapsto& \Ba ^{*}T^{*}(\tau-\cdot)z.
\end{eqnarray*}
\item The Gramian controllability operator is
\begin{equation*}
Q_{\tau \delta*} = G_{\tau\delta}G_{\tau\delta}^{*}= \int_{\tau-\delta}^{\tau}T(\tau-t)\Ba \Ba ^{*}T^{*}(\tau-t)dt.
\end{equation*}
\end{enumerate}

In general, for linear bounded operator $G$ between Hilbert spaces $\mathcal{W}$ and $\mathcal{Z}$, the following lemma holds (see \cite{Bashirov-Kerimov:1997aa,Bashirov-Mahmudov:1999aa, Leiva-Merentes:2013aa}).
\begin{lemma}
The approximate controllability of the linear system \eqref{eq:linear}  on $[\tau-\delta,\tau]$ is equivalent to  any of the following statements
\begin{enumerate}[(a)]
\item $\overline{Rang(G_{\tau\delta})}=\Z^{1}.$
\item $\ker(G_{\tau\delta}^{*})={0}.$
\item For $0\neq z \in\ Z^{1},  \  \ \inner{ Q_{\tau\delta}z,z}>0.$
\end{enumerate}
\end{lemma}

 The controllability of the linear system \eqref{eq:linear} on $[0,\tau]$ was proved  by A. Carrasco and H. Leiva in \cite{Carrasco-Leiva:2013aa}.  Theorem \ref{A1.5}   and Lemma \ref{lema} characterized the controllability of the system \eqref{eq:linear}, their proofs and details can be found in \cite{Bashirov-Kerimov:1997aa,Bashirov-Mahmudov:1999aa, Curtain-Pritchard:1978aa, Curtain-Zwart:1995aa, Leiva-Merentes:2013aa}

 \begin{teo}\label{A1.5} The system \eqref{eq:linear} is approximately
controllable on $[0,\tau]$ if and only if any one of the following conditions hold:
\begin{enumerate}
\item $\ds\lim_{\alpha \to 0^+} \alpha(\alpha I +Q_{\tau\delta}^{*})^{-1}z =0 $.\\
\item If $z\in Z^{1}$,  $0<\alpha \leq 1$ and $u_{\alpha}=G_{\tau\delta}^{*}(\alpha I +
Q_{\tau\delta}^{*})^{-1}z$, then
$$
G_{\tau\delta}u_{\alpha}=z - \alpha(\alpha I+ Q_{\tau\delta})^{-1}z \quad \mbox{and} \quad
\displaystyle\lim_{\alpha\to 0}G_{\tau\delta}u_{\alpha}=z.
$$
Moreover, for each $v\in L^{2}([\tau-\delta,\tau];\U)$, the sequence of controls
$$
u_{\alpha}=G_{\tau\delta}^{*}(\alpha I +
Q_{\tau\delta}^{*})^{-1}z + (v-G_{\tau\delta}^{*}(\alpha I +
Q_{\tau\delta}^{*})^{-1}G_{\tau\delta}v),
$$
satisfies
$$
G_{\tau\delta}u_{\alpha}=z-\alpha(\alpha I +
Q_{\tau\delta}^{*})^{-1}(z-G_{\tau\delta}v)
\quad
\mbox{and}
\quad
\displaystyle\lim_{\alpha\to 0}G_{\tau\delta}u_{\alpha}=z,
$$
with the error $E_{\tau\delta}z=\alpha(\alpha I +
Q_{\tau\delta})^{-1}(z+G_{\tau\delta}v),\;\alpha\in(0,1].
$
\end{enumerate}
\end{teo}
Theorem \ref{A1.5} indicates that the family of linear operators $
\Gamma_{\tau\delta}=G_{\tau\delta}^{*}(\alpha I +
Q_{\tau\delta}^{*})^{-1}
$
is an approximate right inverse for the $G_{\tau\delta}$, in
the sense that
$$
\displaystyle\lim_{\alpha\longrightarrow 0}G_{\tau\delta}\Gamma_{\tau\delta}=I,
$$
in the strong topology.

\begin{lemma}\label{lema}
$Q_{\tau\delta}> 0$, if and only if, the linear system \eqref{eq:linear} is controllable on $[\tau-\delta, \tau]$.
Moreover, for  given initial state $y_0$ and  final state $z_{1}$, there exists a sequence of controls $\{u_{\alpha}^{\delta}\}_{0 <\alpha \leq 1}$ in the space $L^2(\tau-\delta,\tau;\U)$, defined by
$$
u_{\alpha}=u_{\alpha}^{\delta}= G_{\tau\delta}^{*}(\alpha I+ G_{\tau\delta}G_{\tau\delta}^{*})^{-1}(z_{1} - T(\tau)y_0),
$$
such that the solutions $y(t)=y(t,\tau-\delta, y_0, u_{\alpha}^{\delta})$ of the initial value problem
\begin{equation}\label{IVL}
\left\{
\begin{array}{l}
y'=\Aa y+\Ba u_{\alpha}(t), \ \  y \in \Z^{1}, \ \ t>0,\\
y(\tau-\delta) = y_0,
\end{array}
\right.
\end{equation}
satisfies
\begin{equation}\label{eq:limit}
\lim_{\alpha \to 0^{+}}y(\tau)
= \lim_{\alpha \to 0^{+}}\left(T(\delta)y_0 + \int_{\tau-\delta}^{\tau}T(\tau-s)\Ba u_{\alpha}(s)ds \right)= z_{1}.
\end{equation}
\end{lemma}

\section{Controllability of the Semilinear System}\label{sec:semilineal}

This section is devoted to prove the main result of this paper,   the approximate controllability of the beam equation  (Theorem \ref{main}), which is  it is equivalent to prove the controllability of the abstract system \eqref{eq:abstract}  under the condition \eqref{eq:bound}. Recall 

\begin{definition} \label{def2}{\rm (}{\sf Approximate Controllability}{\rm )}  The system \eqref{eq:abstract} is said
to be approximately controllable on $[0,\tau]$
if  for every $\epsilon>0$,   every
$\Phi\in \C\left(-r,0;  \Z^{1} \right)$ and a given initial state $z_{1}\in \Z^{1}$ there exists $u\in L^{2}(0,\tau;\U)$, such that,  the corresponding  mild solution
\begin{align}\label{eq:mild}
z^{u}(t) = & \ds T(t)\Phi(0)+\int_{0}^{t}T(t-s)\left[\Ba u(s)+\left(\int_{0}^{s}\Ga(s,l,z(l-r))dl\right)\right]ds \\
& \;+   \ds \int_{0}^{t}T(t-s)\Fa(s,z(s-r),u(s))ds  +   \sum_{0 < t_k < t} T(t-t_k )\Ia_{k}(t_k,z(t_k), u(t_k)),  \nonumber
\end{align}
satisfies $z(0)=\Phi(0)$ and
\begin{equation}\label{eq:goal}
\norm{ z^u(\tau) - z_{1}}_{\Z^1}<\epsilon.
\end{equation}

\end{definition}

The approach to obtain \eqref{eq:goal} consist in construct a sequence of controls conducting the system from the initial condition $\Phi$ to a small ball around $z_1.$ This is achieved taking advantage of the delay, which  allows us to pullback the corresponding family of  solutions to a fixed trajectory in short time interval. Now, we are ready to present the proof of our main result

\begin{teo} \label{main}
Under the condition  \eqref{eq:f1} the impulsive semilinear beam equation with memory and delay \eqref{eq:beam}-\eqref{eq:initial} is approximately
controllable on $[0,\tau]$.
\end{teo}
\dem Let $\epsilon>0$, and given $\Phi\in \mathcal{C}$ and a final state $z_{1}$. By section \ref{sec:formulation}, we have that the semilinear beam equation in consideration can be represented as the abstract system \eqref{eq:abstract}  under the condition \eqref{eq:bound}.  Thus,
consider any $u\in L^{2}([0,\tau];\U)$ and the corresponding mild solution \eqref{eq:mild}  of the initial value problem \eqref{eq:abstract}, denoted by $z(t)=z(t,0,\Phi,u)$. 

\noindent For $0\leq\alpha \leq 1,$ define the control $u_{\alpha}^{\delta}\in L^{2}([0,\tau];\U)$  as follows
$$
u_{\alpha}^{\delta}(t)=\left\{\begin{array}{ccl}
                         u(t), &&0\leq t\leq \tau-\delta, \\
                         u_{\alpha}(t), &\quad& \tau-\delta\leq t\leq \tau,
                       \end{array}\right.
$$
with $
u_{\alpha}= \Ba^{*}T^{*}(\tau-t)(\alpha I+ G_{\tau\delta}G_{\tau\delta}^{*})^{-1}(z_{1} - T(\delta)z(\tau-\delta)).
$
For,
 $0<\delta<\tau-t_{p}$ its corresponding mild solution at time $\tau$ can be written as follows:
 \begin{eqnarray*}
 \ds
z^{\delta,\alpha}(\tau) &=&
 \ds T(\tau)\Phi(0) +\int_{0}^{\tau}T(\tau-s)
	\left[ \Ba   u_{\alpha}^{\delta} (s)
		+ \int_0^s \Ga(z^{\delta,\alpha}(l-r))dl\right]ds+ \\
	&&+  \int_{0}^{\tau}T(\tau-s)\Fa(s,z^{\delta,\alpha}(s-r),u_{\alpha}^{\delta}(s))ds+
\sum_{0 < t_k < \tau} T(t-t_k )\Ia_{k}(t_k,z^{\delta,\alpha}(t_k), u_{\alpha}^{\delta}(t_k))\\
&=&T(\delta)\left\{T(\tau-\delta)\Phi(0)
+\int_{0}^{\tau-\delta}T(\tau-\delta-s) \left(\Ba   u_{\alpha}^{\delta} (s)+\Fa(s,z^{\delta,\alpha}(s-r),u_{\alpha}^{\delta}(s))\right)ds\right.\\
&&\qquad\quad+\int_{0}^{\tau-\delta}T(\tau-\delta-s)  \int_0^s \Ga(s,l, z^{\delta,\alpha}(l-r))dlds\\
&&\qquad\quad\left.+ \sum_{0 < t_k < \tau-\delta} T(t-\delta-t_k )\Ia_{k}(t_k,z^{\delta,\alpha}(t_k), u_{\alpha}^{\delta}(t_k))\right\}+\\
&& + \int_{\tau-\delta}^{\tau}T(\tau-s)\left(\Ba u_{\alpha}(s)+
\Fa(s,z^{\delta,\alpha}(s-r),u_{\alpha}^{\delta}(s))+\int_0^s\Ga(s,l,z^{\delta,\alpha}(l-r))dl\right)ds.
\end{eqnarray*}
Therefore,
\begin{align*}
z^{\delta,\alpha}(\tau)  = &T(\delta)z(\tau-\delta)+ \int_{\tau-\delta}^{\tau}T(\tau-s)\left(\Ba u_{\alpha}(s)+
\Fa(s,z^{\delta,\alpha}(s-r),u_{\alpha}^{\delta}(s))\right)ds \\
&+ \int_{\tau-\delta}^{\tau}T(\tau-s)\int_0^s\Ga(s,l,z^{\delta,\alpha}(l-r))dlds.
\end{align*}
Observing that the corresponding solution $y^{\delta,\alpha}(t)=y(t,\tau-\delta,z(\tau-\delta),u_{\alpha})$ of the initial value problem \eqref{IVL} at time $\tau$ is:
$$
y^{\delta,\alpha}(\tau)=T(\delta)z(\tau-\delta)+ \int_{\tau-\delta}^{\tau}T(\tau-s)\Ba_{\varpi} u_{\alpha}(s)ds,
$$
yields,
$$
z^{\delta,\alpha}(\tau)-y^{\delta,\alpha}(\tau)=
\int_{\tau-\delta}^{\tau}T(\tau-s)\left(\int_{0}^{s}\Fa(s,z^{\delta,\alpha}(s-r),u_{\alpha}^{\delta}(s))+\Ga(s,l,z^{\delta,\alpha}(l-r))dl)\right)ds,
$$
 and together with condition \eqref{eq:bound}, we obtain 
\begin{align*}
  \norm{ z^{\delta,\alpha}(\tau)-y^{\delta,\alpha}(\tau)} & \leq \int_{\tau-\delta}^{\tau} \norm{ T(\tau-s)}\left( \tilde{a}\norm{\Phi(s-r)}+\tilde{b}\right)ds \\
   & +  \int_{\tau-\delta}^{\tau}\norm{ T(\tau-s)}\int_{0}^{s}\norm{\Ga(s,l,z^{\delta,\alpha}(l-r))}dlds.
\end{align*}
Observe that
 $0< \delta< r$ and $\tau-\delta \leq s\leq \tau$, thus $$l-r \leq s-r \leq \tau-r< \tau-\delta.
 $$
Therefore,
$
z^{\delta,\alpha}(l-r)=z(l-r) $  and $ z^{\delta,\alpha}(s-r)=z(s-r),
$ implying that for $\epsilon>0,$ there exists $\delta>0$ such that
\begin{align*}
 \norm{z^{\delta,\alpha}(\tau)-y^{\delta,\alpha}(\tau)} & \leq \int_{\tau-\delta}^{\tau}\norm{ T(\tau-s)}\left( \tilde{a}\norm{z(s-r)}+\tilde{b}\right)ds \\
  &\quad +  \int_{\tau-\delta}^{\tau}\norm{T(\tau-s)}\int_{0}^{s}\norm{ \Ga(s,l,z(l-r))} dlds  \\
   & <  \displaystyle\frac{\epsilon}{2}.
\end{align*}
Additionally, for $0<\alpha <1$, Lemma \ref{lema}  \eqref{eq:limit} yields
$$
 \norm{ y^{\delta,\alpha}(\tau)-z_{1}}  <  \frac{\epsilon}{2}.
$$
Thus,
$$
\begin{array}{lll}
 \norm{ z^{\delta,\alpha}(\tau)-z_{1}}  & \leq &  \norm{ z^{\delta,\alpha}(\tau)-y^{\delta,\alpha}(\tau)} +  \norm{ y^{\delta,\alpha}(\tau)-z_{1}}  <  \frac{\epsilon}{2}+ \frac{\epsilon}{2}=\epsilon,
\end{array}
$$
which completes our proof.
\section{Final Remarks}\label{final}
\noindent

We believe this technique  can be applied for  controlling diffusion processes systems involving compact semigroups. In particular, our result can be formulated in a more general setting  for the semilinear evolution equation with impulses, delay and memory  in a Hilbert space $\Z$
\begin{equation*}\label{eq:class}
\left\{
\begin{array}{lr}
z' =  -\Aa z+ \Ba u + \ds \int_{0}^{t}\Ga(t,s ,z(s-r))ds +  \Fa(t,z_{t}(-r),u(s)) ,& z\in Z^{1},\; t\geq 0, \\
z(s)  = \Phi(s), & s \in [-r,0],\\
z(t_{k}^{+})  = z(t_{k}^{-})+\Ia_{k}(t_k, z(t_{k}),u(t_{k})), &  k=1,2, \dots, p,
\end{array}
\right.
\end{equation*}
where $u\in L^{2}(0,\tau;\U)$, $\U$ is another Hilbert space, $\Ba :\U \longrightarrow \Z$ is a bounded linear operator, $\Ia_{k}, \Fa:[0, \tau]\times \C(-r,0; \Z) \times \U \rightarrow \Z$, $\Aa :D(\Aa) \subset \Z \rightarrow \Z$ is an unbounded linear operator in $\Z$ that generates a strongly continuous semigroup \cite[Lemma 2.1]{Leiva:2003aa}
\begin{equation*}\label{damp2}
  T(t)z =\sum_{nj=1}^{\infty}e^{\Aa_{j}t}P_jz%
  \mbox{, } \ \ z\in \Z \mbox{, } \ \ t \geq 0,     
\end{equation*}
where  $\left\{ P_j\right\} _{j \geq 0}$ is a complete family of orthogonal projections in the Hilbert space $\Z$ and
\begin{equation*}
\|\Fa(t,\Phi,u) \|_{\Z}  \leq   \tilde{a} \|\Phi(-r)\|_{\Z} +\tilde{b}, 
\end{equation*}
for all $(t, \Phi, u) \in [0, \tau]\times \C(-r,0;  \Z ) \times \U.$

\begin{center} {\sc Acknowledgments}
\end{center}
The authors  are  thankful to the anonymous referees for valuable comments 
that help  improve the quality of the paper.  This work has been supported by Louisiana State University, Universidad YachayTech  and  Universidad Centroccidental Lisandro Alvarado.

\end{document}